\documentclass[runningheads]{llncs}
\usepackage[T1]{fontenc}
\usepackage{graphicx}
\usepackage{amsmath}
\usepackage{empheq}
\usepackage{bbm}
\usepackage{amssymb}
\usepackage{mathrsfs}  
\usepackage{bm}
\usepackage{circuitikz}
\newcommand{\giantzero}{\mbox{\normalfont\Huge0}}

\begin{document}
\title{Expressiveness and Structure Preservation in Learning  Port-Hamiltonian Systems}
%
%
\author{Juan-Pablo~Ortega and Daiying Yin}
%
%
\institute{Division of Mathematical Sciences, Nanyang Technological University, Singapore\\
\email{\texttt{Juan-Pablo.Ortega@ntu.edu.sg}} and {\texttt{yind0004@e.ntu.edu.sg}}}
\titlerunning{Learning Port-Hamiltonian Systems}
\maketitle              
\begin{abstract}
	A well-specified parametrization for single-input/single-output (SISO) linear port-Hamiltonian systems amenable to structure-preserving supervised learning is provided. The construction is based on controllable and observable normal form Hamiltonian representations for those systems, which reveal fundamental relationships between classical notions in control theory and crucial properties in the machine learning context, like structure-preservation and expressive power. The results in the paper suggest parametrizations of the estimation problem associated with these systems that amount, at least in the canonical case, to unique identification and prove that the parameter complexity necessary for the replication of the dynamics is only $\mathcal{O}(n)$ and not $\mathcal{O}(n^2)$, as suggested by the standard parametrization of these systems. \\
	
\end{abstract}

\section{Introduction}
%
%
%
%
Machine learning has experienced substantial development in recent years due to significant advances in algorithmics and a fast growth in computational power. In physics and engineering, machine learning is called to play an essential role in predicting and integrating the equations associated with physical dynamical systems. Nevertheless, for physics-related problems, like in mechanics or optics, it is natural to build into the learning algorithm any prior knowledge that we may have about the system based on physics' first principles. This may include specific forms of the laws of motion, conservation laws, symmetry invariance, as well as other underlying geometric and variational structures. This observation regarding the construction of structure-preserving schemes, has been profusely exploited with much success before the emergence of machine learning in the field of numerical integration \cite{gonzalez2000time}, \cite{marsden_west_2001}, \cite{LeimkuhlerBook}, \cite{mclachlan2006geometric}. Many examples in that context show how the failure to maintain specific conservation laws can lead to physically inconsistent solutions. The translation of this idea to the context of machine learning has lead to the emergence of a new domain collectively known as {\it physics-informed machine learning} (see \cite{raissi2017physics}, \cite{wu2018physics}, \cite{karniadakis2021physics} and references therein).

The approach that we propose in this contribution differs from the references mentioned above in two ways. First, these methods are designed to learn the state evolution of Hamiltonian systems, whereas our approach focuses on {\it learning the input-output dynamics of port-Hamiltonian systems while forgetting about the physical state space}. As will be introduced later, these systems have an underlying Dirac structure that describes the geometry of numerous physical systems with external inputs \cite{vanderSchaft2014} and includes the dynamics of the observations of Hamiltonian systems as a particular case. This is a significant difference with respect to the results available in the literature, which mostly deal with autonomous Hamiltonian systems on which one assumes access to the entire phase space and not only to its observations. Second, instead of a general nonlinear system for which only approximation error can be possibly estimated, we consider, as a first approach {\it exclusively linear systems}, for which we can obtain explicit representations in normal form which allows us to propose a structure-preserving learning paradigm with a provable minimal parameter space (in the canonical case).

 The contributions in this paper are mostly based on classical techniques in control theory, the Cayley-Hamilton theorem, and symplectic linear algebra that allow us to define the notion of {\it normal-form controllable and observable Hamiltonian representations}. Our main result (Theorem \ref{main_theorem}) shows the existence of system morphisms that allow us to represent any linear port-Hamiltonian system in normal form as the image of a normal-form controllable Hamiltonian representation of the same dimension. An obvious observation is that since the original port-Hamiltonian system and the new linear system are linked by a system morphism, the image of the input/output relations of the latter is the input/output relations of the former. In particular, the new system can be used to learn to reproduce the input/output dynamics of the original port-Hamiltonian system (for a subspace of initial conditions) and {\it this learning paradigm is structure-preserving by construction}. Analogously, in the same result, we introduce another type of system morphisms that link any linear port-Hamiltonian system to some normal-form observable Hamiltonian representation of the same dimension. Consequently, the input-ouput relations of the original port-Hamiltonian system with respect to any initial condition can be reproduced using the observable Hamiltonian representation. 
 
The controllable and observable representations are closely related to each other, and both system morphisms become isomorphisms for canonical port-Hamiltonian systems. However, for the purpose of learning a general port-Hamiltonian system that may not be canonical, we reveal that there is a trade-off between the {\it structure-preserving property} and the {\it expressive power}. These results establish a strong link between classical notions in control theory, e.g. controllability and observability, and those in machine learning, e.g. structure-preservation and expressive power.

The results in the paper point at parametrizations of the estimation problem associated with these systems that amount, at least in the canonical case, to unique identification and prove that the parameter complexity necessary for the replication of the dynamics is only $\mathcal{O}(n)$ and not $\mathcal{O}(n^2)$, as suggested by the standard parametrization of these systems. We shall sketch some results that point to a characterization of the unique identification problem in terms of groupoid/group orbit spaces that will be studied in detail in a forthcoming publication.

\section{Preliminaries}
\label{preliminary}

\subsubsection{State-space systems and morphisms.} 
A continuous time state-space system is given by the following two equations 
\begin{equation}
	\label{state space system}
	\left\{
	\begin{aligned}
		\dot{{\bf z}}&=F({\bf z}, u),\\
		y&=h({\bf z}),
	\end{aligned} \right.
\end{equation} 
where $u \in \mathcal{U}$ is the {\it input}, ${\bf z} \in \mathcal{Z}$ is the {\it internal state},  $y \in {\cal Y} $ is the {\it output}, $F:\mathcal{Z}\times\mathcal{U}\rightarrow\mathcal{Z}$ is called the {\it state map}, and $h : {\cal Z} \rightarrow {\cal Y} $ is the {\it readout} or {\it output map}. The first equation is called the {\it state equation} while the second is usually called the {\it observation equation}. State-space systems will sometimes be denoted using the triplet $(\mathcal{Z}, F, h)$. The solutions of \eqref{state space system} (when available and unique) yield an {\it input/output map} or {\it filter} $U_{F,h}: C^{1}(I, \mathcal{U}) \times {\cal Z} \longrightarrow C^{1}(I, {\cal Y}) $ that is by construction causal and time-invariant. The filter $U_{F,h} $ associates to each pair $(\mathbf{u}, z _0) \in C^{1}(I, \mathcal{U}) \times {\cal Z}  $ the solution $U_{F,h}(\mathbf{u}, z _0) \in C^{1}(I, {\cal Y})$ of the non-autonomous differential equation \eqref{state space system}. Different state-space systems may induce the same input/output system. This is what we call the {\it identification problem}, which is related to the next definition (see \cite{RC16}).

\begin{definition}
	\label{morphism}
	A map $f:\mathcal{Z}_1\rightarrow\mathcal{Z}_2$ is called a {\it system morphism}  between the continuous-time state-space systems $(\mathcal{Z}_1, F_1,h_1)$ and $(\mathcal{Z}_2, F_2,h_2)$ if it satisfies the following two properties:
	\begin{description}
		\item[(i)] {\it System equivariance}: $f(F_1({\bf z}_1,u))=F_2(f({\bf z}_1),u)$, for all ${\bf z}_1\in\mathcal{Z}_1$ and $u\in\mathcal{U}$.
		\item[(ii)] {\it Readout invariance}: $h_1({\bf z}_1)=h_2(f({\bf z}_1))$ for all ${\bf z}_1\in\mathcal{Z}_1$.
	\end{description}
\end{definition}
An elementary but very important fact is that if
$f:\mathcal{Z}_1\rightarrow\mathcal{Z}_2$ is a linear system-equivariant map between $(\mathcal{Z}_1, F_1,h_1)$ and $(\mathcal{Z}_2, F_2,h_2)$ ($\mathcal{Z}_1 $ and $\mathcal{Z}_2 $ are in this case vector spaces) then, for any solution
${\bf z}_1\in C ^1(I, \mathcal{Z}_1)$ of the state equation associated to $F_1 $ with initial condition ${\bf z} _1 ^0 \in {\cal Z}  $ and to the input $\mathbf{u} \in C ^1(I, \mathcal{U})$, its image $f \circ {\bf z} _1 \in C ^1(I, \mathcal{Z}_2)$  is a solution for the state space system associated to $F_2$ with the same input and initial condition $f({\bf z} _1 ^0) \in {\cal Z} _2$. This implies that 
\begin{equation}
\label{congfil}
U_{F _2, h _2}\circ \left(\mathbb{I}_{C ^1(I, \mathcal{U})} \times f \right)=U_{F _1, h _1}.
\end{equation}
This fact has as an important consequence that, in general, input/output systems {\it are not uniquely identified} since all the system-isomorphic state-space systems via a map $f$ yield the same input/output map, with initial conditions related by the map $f$ according to \eqref{congfil}. 

\subsubsection{Hamiltonian and port-Hamiltonian systems.}

The {\it Hamiltonian system} determined by the {\it Hamiltonian function} $H \in C^1(\mathbb{R}^{2n} ) $ is given by the differential equation
\begin{equation}
	\label{hamiltonian_def}
	\dot{{\bf z}}=\mathbb{J}\frac{\partial H}{\partial {\bf z}},
\end{equation}
where  ${\displaystyle \mathbb{J}=
	\begin{bmatrix}
		0 & \mathbb{I}_n\\
		-\mathbb{I}_n & 0
\end{bmatrix}}$ is the so-called the {\it canonical symplectic matrix}. 
A {\it linear} Hamiltonian system is determined by a quadratic Hamiltonian function $H({\bf z})=\frac{1}{2}{\bf z}^{T}Q{\bf z}$, where ${\bf z}\in\mathbb{R}^{2n}$ and $Q \in \mathbb{M}_{2n}$ is a square matrix that without loss of generality can be assumed to be symmetric. In this case, Hamilton's equations (\ref{hamiltonian_def}) reduce to $\dot{{\bf z}}=\mathbb{J}Q{\bf z}$.

\noindent {\it Port-Hamiltonian systems} (see \cite{vanderSchaft2014}) are state-space systems that generalize autonomous Hamiltonian systems to the case in which external signals or inputs control in a time-varying way the dynamical behavior of the Hamiltonian system.  The family of input-state-output port-Hamiltonian systems have the explicit form (see \cite{vanderSchaft2014}):
\begin{equation}
	\label{general ph system}
	\left\{
	\begin{aligned}
		\dot{{\bf x}}&=[J({\bf x})-R({\bf x})]\frac{\partial H}{\partial {\bf x}}({\bf x})+g({\bf x})u,\\
		y&=g^{T}({\bf x})\frac{\partial H}{\partial {\bf x}}({\bf x}),
	\end{aligned} \right.
\end{equation}
where $(u,y)$ is the input-output pair (corresponding to the control and output conjugated ports), $J({\bf x})$ is a skew-symmetric interconnection structure and $R({\bf x})$ is a symmetric positive-definite dissipation matrix. 
Our work concerns {\it linear} port-Hamiltonian systems in  {\it normal form}, that is, the skew-symmetric matrix $J$ is constant and equal to the canonical symplectic matrix $\mathbb{J}$, the Hamiltonian matrix $Q$ is symmetric positive-definite, and the energy dissipation matrix $R=0$, in which case \eqref{general ph system} takes the form:
\begin{equation}
	\label{ph_definition}
	\left\{
	\begin{aligned}
		\dot{{\bf z}}&=\mathbb{J}Q{\bf z}+\mathbf{B}u,\\
		y&=\mathbf{B}^{T}Q{\bf z},
	\end{aligned} \right.
\end{equation} 
with ${\bf z}\in\mathbb{R}^{2n}$, $u,y\in\mathbb{R}$, and where $\mathbf{B}\in\mathbb{R}^{2n}$ specifies the interconnection structure simultaneously at the input and output levels. In all that follows, we denote the family of state-space systems that can be written as \eqref{ph_definition} with the symbol ${PH}_n$, and we refer to the elements of this set as {\it normal form port-Hamiltonian systems}. All these systems have the existence and uniqueness of solutions property and hence determine a family of input/output systems that will be denoted by $\mathcal{PH}_n$. By definition, the systems in ${PH}_n$ are fully determined by the pairs $(Q,\mathbf{B})$, and hence we define the parameter space 
\begin{equation*}
\Theta_{{PH}_n}:=\left\{(Q,\mathbf{B})|0<Q\in\mathbb{M}_{2n}, Q=Q^T, \mathbf{B}\in\mathbb{R}^{2n}\right\}.
\end{equation*}
Let $\theta_{PH_n}:\Theta_{PH_n} \rightarrow PH_n $ the map that associates to the parameter $(Q,\mathbf{B}) \in \Theta_{PH_n} $ the corresponding port-Hamiltonian state space system \eqref{ph_definition}.
We emphasize that two different elements in $\Theta_{{PH}_n} $ can determine the same input/output system in $\mathcal{PH}_n$; said differently, the parameter set $\Theta_{{PH}_n}$ does not uniquely specify the elements in $\mathcal{PH}_n$.

\subsubsection{Williamson's normal form.}
\label{Williamson normal form section}
The following classical result can be found in \cite{williamson1936algebraic}, \cite{williamson1937normal}, \cite{Ikramov2018}, \cite{de2006symplectic}.
\begin{theorem}
	\label{Williamson's normal form}
	Let $M \in \mathbb{M}_{2n}$ be a positive-definite symmetric real matrix. Then
	\begin{description}
		\item[(i)] There exists a symplectic matrix $S\in Sp(2n, \mathbb{R})$ such that ${\displaystyle M=S^{T}\begin{bmatrix}
				D&0\\
				0&D
		\end{bmatrix}}S$, with $D= {\rm diag}({\bf d})$ a $n$-dimensional diagonal matrix with positive entries and ${\bf d}=\left(d_1, \ldots, d_n\right)^T$
		.\\
		\item[(ii)] The values $d_1, \ldots, d_n$ are independent, up to reordering, on the choice of the symplectic matrix $S$ used to diagonalize $M$.
	\end{description}
\end{theorem}

We use the notation $D={\rm diag}({\bf d})$ to denote that $D$ is a diagonal matrix with diagonal entries given by the vector ${\bf d}=(d_1,\dots,d_n)^T$. The elements $d_i$ in the above theorem are called the {\it symplectic eigenvalues} of $M$.

\section{Controllable and Observable Hamiltonian Systems}

\begin{definition}\label{definition_of_controllable}
	Given ${\bf d}=\left(d_1, \ldots, d_n\right)^T\in {\Bbb R}^n$, with $d_i>0$, and ${\bf v}\in\mathbb{R}^{2n}$, we say that a $2n$-dimensional linear state space system is a {\it controllable Hamiltonian} (respectively, {\it observable Hamiltonian}) representation if it takes the form
	\begin{equation}
		\label{proposed_system_learner}
		\left\{
		\begin{aligned}
			\dot{{\bf s}}&=g^{ctr}_1({\bf d})\cdot {\bf s}+\left(
			0 , 0 ,\cdots, 0, 1\right)^{T}\cdot u,\\
			y&=g^{ctr}_2({\bf d}, {\bf v})\cdot {\bf s},
		\end{aligned} \right. 
\left(\mbox{resp., } \left\{
	\begin{aligned}
		\dot{{\bf s}}&=g^{obs}_1({\bf d})\cdot {\bf s}+g^{obs}_2({\bf d},{\bf v})\cdot u,\\
		y&=\left(
		0 , 0 ,\cdots, 0, 1\right)\cdot {\bf s},
	\end{aligned} \right. \right)
	\end{equation} 
	where $g^{ctr}_1({\bf d}) \in \mathbb{M}_{2n }$ and $g^{ctr}_2({\bf d}, {\bf v}) \in \mathbb{M}_{1,2n}$ (respectively, $g^{obs}_1({\bf d}) \in \mathbb{M}_{2n }$ and $g^{obs}_2({\bf d}, {\bf v}) \in \mathbb{R}^{2n}$) are constructed as follows:
	\begin{description}
		\item [(i)] Given ${\bf d} \in {\Bbb R}^n$, let  $\left\{a _0, a _1, \ldots , a_{2n-1}\right\}$ be the real coefficients that make ${\lambda}^{2n}+\sum_{i=0}^{2n-1}a_i\cdot{\lambda}^{i}=({\lambda}^2+d_1^2)({\lambda}^2+d_2^2)\dots({\lambda}^2+d_n^2)$ an equality between the two polynomials in $\lambda$. Let $a_{2n}=1$ by convention. Note that the entries $a_{i}$ with an odd index $i $ are zero. Define: \begin{equation*}
			g^{ctr}_1({\bf d}):=\begin{bmatrix}
				0&1&0&\dots&0 \\
				0&0&1&\dots&0 \\
				\vdots&\vdots&\ddots&\vdots&\vdots\\
				0&0&0&\dots&1 \\
				-a_0&-a_1&-a_2&\dots&-a_{2n-1}
			\end{bmatrix}_{2n\times 2n},
		\end{equation*}
(respectively, $g^{obs}_1({\bf d})=g^{ctr}_1({\bf d}) ^{\top}$).
		\item[(ii)] Given ${\bf d}$ and ${\bf v}$, then \begin{equation*}
			\begin{aligned}
				g^{ctr}_2({\bf d},{\bf v})&:=\begin{bmatrix}
					0\:\:c_{2n-1}\:\:0\:c_{2n-3}\:\hdots\:\:0\:\:c_{1}
				\end{bmatrix}, \quad \mbox{(resp., $g^{obs}_2({\bf d},{\bf v})= g^{ctr}_2({\bf d},{\bf v}) ^{\top}$)} 
			\end{aligned}
		\end{equation*} where 
		\begin{equation*}
			c_{2k+1}={\bf v}^{T}\begin{bmatrix}
				F_k&0\\
				0&F_k
			\end{bmatrix} {\bf v},
		\end{equation*} for $k=0,\dots,n-1$, and
		\begin{equation*}
			F_k=\begin{bmatrix}
				f_1&&\\
				&f_2&&\giantzero\\
				&&\ddots\\
				&\giantzero&&f_{n-1}\\
				&&&&f_n\\
			\end{bmatrix}
		\end{equation*}
		with $f_l=d_l\cdot\sum_{\substack{j_1,\dots,j_k\neq l\\1\leq j_1<\dots<j_k\leq n}}\big(d_{j_1}d_{j_2}\cdots d_{j_k}\big)^2$, $l=1,\dots,n$.
	\end{description}
	We denote $CH_n$ (respectively, $OH_n$) the set of all systems of the form (\ref{proposed_system_learner}), and we call them {\it controllable Hamiltonian} (respectively, {\it observable Hamiltonian}) representations. The symbol $\mathcal{CH}_n$ (respectively, $\mathcal{OH}_n$) denotes the set of input/output systems induced by the state space systems in $CH_n$ (respectively, $OH_n$). We emphasize that the elements of both $CH_n$ and $OH_n$ can be parameterized with the set 
	$$\Theta_{CH_n}=\Theta_{OH_n}:=\left\{({\bf d},{\bf v})|d_i>0, {\bf v}\in\mathbb{R}^{2n}\right\}.$$
\end{definition}
The maps $\theta_{CH_n}:\Theta_{CH_n} \rightarrow CH_n $ and $\theta_{OH_n}:\Theta_{OH_n} \rightarrow OH_n $ associate to each parameter the corresponding state space system. Note that the elements in $CH _n $ (respectively, in $OH _n $) of the form (\ref{proposed_system_learner}) are in canonical controllable (respectively, observable) form in the sense of \cite{sontag:book} and they are hence controllable (respectively, observable). Our main result establishes a relationship between port-Hamiltonian systems and controllable (respectively, observable) Hamiltonian representations as defined above, which will be used later on for considerations on the structure preservation and expressiveness modelling of $PH _n $.

\begin{theorem}
\label{main_theorem}
\begin{description}
\item [(i)] There exists, for each $S\in Sp(2n,\mathbb{R})$,  a map
\begin{equation*}
	\begin{array}{cccc}
		\varphi_S: &CH _n& \longrightarrow & {PH}_n\\
		&\theta_{CH _n}({\bf d},{\bf v}) &\longmapsto & \theta_{PH_n}\left(S^{T}\begin{bmatrix}
			D&0\\
			0&D
		\end{bmatrix}S, S^{-1}{\bf v}\right),
	\end{array}
\end{equation*}
with $D={\rm diag}({\bf d})$, such that the controllable Hamiltonian system $\theta_{CH _n}({\bf d},{\bf v}) \in CH_n $  and the port-Hamiltonian image $\varphi_S\left(\theta_{CH _n}({\bf d},{\bf v})\right) \in {PH}_n$ are linked by a linear system morphism $f_S^{({\bf d},{\bf v})}:\mathbb{R}^{2n}\rightarrow\mathbb{R}^{2n}$.
\item [(ii)] Given a port-Hamiltonian system $\theta_{PH _n}(Q,\mathbf{B})\in {PH}_n$,  there exists an explicit linear system morphism $f^{(Q,\mathbf{B})}:\mathbb{R}^{2n}\rightarrow\mathbb{R}^{2n}$ between the state space of $\theta_{PH _n}(Q,\mathbf{B})\in {PH}_n$ and that of an observable Hamiltonian system $\theta_{OH _n}({\bf d},{\bf v})\in OH _n$, where $({\bf d},{\bf v}) \in \Theta_{OH _n}$ is determined by the Williamson's normal form decomposition of $Q$ determined by $S\in Sp(2n,\mathbb{R})$, that is,  ${\displaystyle Q=S^{T}\begin{bmatrix}
			D&0\\
			0&D
	\end{bmatrix}}S$, $D={\rm diag}({\bf d})$ and ${\bf v}=S\cdot \mathbf{B}$.
\end{description}
\end{theorem}

\subsubsection{Controllability, observability, and invertibility.}
The linear system morphism $f_S^{({\bf d},{\bf v})}:\mathbb{R}^{2n}\rightarrow\mathbb{R}^{2n}$ in part {\bf (i)} is implemented by a matrix $L$, that is,  ${\bf z}=f_S^{({\bf d},{\bf v})}({\bf s})=L{\bf s}$ and $L$ can be explicitly computed once $({\bf d},{\bf v}) $ has been fixed. Moreover, $L$ can be transformed by elementary column operations into the controllability matrix of $\varphi_S\left(\theta_{CH _n}({\bf d},{\bf v})\right) \in {PH}_n$. Consequently, the condition of $L$ being invertible, i.e. that the two systems in $CH _n  $ and $PH _n $ linked by $\varphi _S$ being isomorphic, is equivalent to the controllability matrix of $ \varphi_S\left(\theta_{CH _n}({\bf d},{\bf v})\right) $ having full rank (regardless of the choice of $S\in Sp(2n,\mathbb{R})$), which is again equivalent to $ \varphi_S\left(\theta_{CH _n}({\bf d},{\bf v})\right) $ being canonical since controllability and observability are intertwined concepts in the linear port-Hamiltonian category. Indeed, it can be proved (see \cite{Medianu2013}) that if a linear port-Hamiltonian system without dissipation is controllable and $\det(Q)\neq0$, then it is also observable. Conversely, if it is observable, then this implies that $\det(Q)\neq0$ and it is also controllable. As it is customary in systems theory, we say a linear port-Hamiltonian system in normal form is {\it canonical} if it is both controllable and observable. In view of the results that we just recalled, the condition $\det(Q)\neq 0$ which is part of the definition of $PH _n$ implies that in this context either controllability or observability is equivalent to the system being canonical. Thus, in the remaining, we only refer to controllability. 

The systems in $CH _n$ are by construction in controllable canonical form and are therefore always controllable. If the image system $ \varphi_S\left(\theta_{CH _n}({\bf d},{\bf v})\right)$ that we want to learn is controllable (or equivalently, observable), then by the previous argument, $L$ is necessarily an invertible matrix which means that $ \theta_{CH _n}({\bf d},{\bf v})$ and $ \varphi_S\left(\theta_{CH _n}({\bf d},{\bf v})\right) $ are necessarily isomorphic systems by construction. As a consequence, $\theta_{CH _n}({\bf d},{\bf v}) $ is in such case not only controllable but also observable. 

\subsubsection{Application to structure-preserving system learning.}
As a corollary of the previous result, we can use controllable Hamiltonian systems to learn port-Hamiltonian systems in an efficient and structure-preserving fashion. Indeed, given a realization of a port-Hamiltonian system, a system of the type $\theta_{CH_n}({\bf d},{\bf v})\in CH_n$ can be estimated using an appropriate loss. A controllable Hamiltonian representation  is more advantageous than the original port-Hamiltonian one for two reasons:
\begin{description}
\item [(i)] The {\it model complexity} of the controllable Hamiltonian representation is only of order $\mathcal{O}(n)$, as opposed to $\mathcal{O}(n^2)$ for the original port-Hamiltonian one.
\item [(ii)] This learning scheme is automatically {\it structure-preserving}. Indeed, once a system $\theta_{CH _n}({\bf d},{\bf v})\in CH_n$ has been estimated for a given realization, we have shown that there exists a family of linear morphisms, each of which is between the state space of $\theta_{CH _n}({\bf d},{\bf v})$ and some $\theta_{PH _n}(Q,{\bf B})\in  {PH}_n$, such that any solution of (\ref{proposed_system_learner}) is automatically a solution of some system in ${PH}_n$. Hence, {\it even in the presence of estimation errors} for the parameters  $({\bf d},{\bf v})\in \Theta_{CH _n}$, its solutions correspond to a port-Hamiltonian system. Hence, this structure is {\it preserved} by the learning scheme. 
\end{description}

\subsubsection{System learning and expressive power.} 
There is an important relation between the controllability of a system in ${PH}_n$ and the expressive power of the corresponding representation in $CH_n$. Indeed, if a system of the form (\ref{ph_definition}) in ${PH}_n$ is controllable (i.e. $f_S^{({\bf d},{\bf v})}$ invertible), its preimage by the map $\varphi_S  $  in  $CH _n $ reproduces all possible solutions of (\ref{ph_definition}), which amounts to the learning scheme based on the systems in $CH _n $ having full expressive power. However, if $\theta_{PH _n}(Q,{\bf B})\in PH_n$ fails to be controllable (i.e. $f _S^{({\bf d},{\bf v})} $ not invertible), then the full expressive power is not guaranteed. As a rule of thumb, the more controllable a system of the type $\theta_{PH _n}(Q,{\bf B})\in \mathcal{PH}_n$ is, the higher the rank of $f _S^{({\bf d},{\bf v})}  $ is, and then the more expressive the corresponding controllable Hamiltonian representation is.

On the other hand, we have shown that for each $\theta_{PH _n}(Q,{\bf B})\in  {PH}_n$, there exists a linear morphism between the state space of $\theta_{PH _n}(Q,{\bf B})\in  {PH}_n$ and some $\theta_{OH _n}({\bf d},{\bf v})$, such that any solution of (\ref{ph_definition}) is automatically a solution of some system in ${OH}_n$. Hence, {\it no matter the original port-Hamiltonian system is controllable or not}, its solution can always be reproduced by an element in $OH_n$ by learning a correct parameter $({\bf d},{\bf v})\in \Theta_{OH_n}$. Hence, this learning scheme has {\it full expressive power}.

\section{Unique identifiability of  Port-Hamiltonian Systems}
\label{quotient}

The main theorem in the previous section provides a structure-preserving learning scheme of optimal complexity. It does not, however, solve the unique identifiability problem for the family of filters $\mathcal{PH} _n $  induced by the linear port-Hamiltonian systems in normal form \eqref{ph_definition} that we have set to learn in this paper. Even though this problem will be studied in detail in a forthcoming publication, we put forward some preliminary results that can be obtained as a consequence of Theorem \ref{main_theorem}.

We recall that the unique identification problem in our case has its source in the fact that system isomorphic elements in $PH _n  $ induce exactly the same input/output system in $\mathcal{PH} _n $. It is easy to show that system isomorphisms induce an equivalence relation on $PH _n  $ that we shall denote by $\sim_{sys}$.

One may hope that  
\begin{equation*}
	 \mathcal{PH} _n \simeq PH _n /\sim_{sys},
\end{equation*}
so that the quotient space $PH _n /\sim_{sys} $ is the space in which unique identifiability is achieved. Unfortunately, this holds only when we are restricting to canonical systems, which form an open and dense subset. More precisely, if we denote by $PH^{can}_n$ the subset of canonical port-Hamiltonian systems in $PH_n$ and denote by $\mathcal{PH}^{can}_n$ the corresponding filters, then it holds that 
\begin{equation*}
	\mathcal{PH}^{can}_n\simeq PH^{can}_n /\sim_{sys}\simeq \Theta^{can}_{{CH} _n}/(S_n\rtimes \mathbb{T}^n).
\end{equation*} Furthermore, this unique identification space can be characterized as the orbit space of certain action of the semi-direct product group $S_n\rtimes \mathbb{T}^n$ ($S_n$ is the permutation group and $\mathbb{T}^n$ is the $n$-tori) on some dense subset $\Theta^{can}_{CH_n}$ of $\Theta_{CH_n}$. The orbit space has a natural smooth structure and global charts that can be used in the numerical implementation of learning problems. 

Generally speaking, in the presence of non-canonical port-Hamiltonian systems, the relationship between $\mathcal{PH}_n$ and $PH_n/\sim_{sys}$ is complicated. In such a scenario, two distinct elements in $PH_n$ that are $\sim_{sys}$-equivalent always induce the same filter in $\mathcal{PH}_n$, whereas a filter in $\mathcal{PH}_n$ could be realized by two elements in $PH_n$ that are not $\sim_{sys}$-equivalent, since a filter identifies the canonical part (i.e. minimal realization) and that part only, see \cite{kalman1963}. Said differently, by going to the quotient space $PH_n/\sim_{sys}$, we remove some redundancies of the set $PH_n$ that yield the same input-output dynamics, but not all. In what follows, we give a characterization of $PH_n/\sim_{sys}$ while we emphasize again that $PH_n/\sim_{sys}$ is not the same as the unique identifiability space $\mathcal{PH}_n$ in general. 

In order to make the quotient space $PH_n/\sim_{sys}$ manageable, it is desirable to characterize it in terms of the parameter set $\Theta_{{PH} _n} $ or, even better, by the less complex $\Theta_{{CH} _n} $. This is indeed possible. More explicitly, an equivalence relation $\sim_{\star} $ on $\Theta_{{CH} _n} $ can be defined such that  
\begin{equation}
	\label{isomorphisms}
	 PH _n /\sim_{sys} \simeq \Theta_{{CH} _n}/\sim_{\star}.
\end{equation}
The equivalence relation $\sim_{\star} $ is defined as: the pairs $({\bf d}_1,{\bf v}_1)$ and $({\bf d}_2,{\bf v}_2)$ in $\Theta_{{CH} _n}$ are $\sim_{\star} $-equivalent, i.e. $({\bf d}_1,{\bf v}_1)\sim_{\star}({\bf d}_2,{\bf v}_2)$, if there exists a permutation matrix $P_{\sigma}\in\mathbb{M}_{n}$ and an invertible matrix $A$ such that, for $D_i={\rm diag}({\bf d}_i)$, $i \in \left\{1,2\right\} $ and $P=\begin{bmatrix}
	P_{\sigma}&0\\0&P_{\sigma}
\end{bmatrix}$, the following conditions hold true:
$$P\begin{bmatrix}
	D_{1}&0\\0&D_{1}
\end{bmatrix}P^{T}=\begin{bmatrix}
	D_{2}&0\\0&D_{2}
\end{bmatrix}, \quad
A^{T}\begin{bmatrix}
	D_{1}&0\\0&D_{1}
\end{bmatrix}A{\bf v}_1=\begin{bmatrix}
	D_{1}&0\\0&D_{1}
\end{bmatrix}{\bf v}_1,$$
$$A\mathbb{J}\begin{bmatrix}
	D_{1}&0\\0&D_{1}
\end{bmatrix}=\mathbb{J}\begin{bmatrix}
	D_{1}&0\\0&D_{1}
\end{bmatrix}A, \quad \mbox{and} \quad
{\bf v}_2=PA{\bf v}_1$$

\subsubsection{Characterization of $PH_n/\sim_{sys}$ in terms of Lie groupoid orbit spaces.} As it is customary, groupoids will be denoted with the symbol $s,t: {\cal G} \rightrightarrows M $ (or simply ${\cal G} \rightrightarrows M $), where $s$ and $t $ are the {\it source} and the {\it target} maps, respectively. Given $m \in M $, the {\it groupoid orbit} that contains this point is given by $\mathcal{O} _m= t \left(s ^{-1}(m)\right) \subset M$. The {\it orbit space} associated to ${\cal G} \rightrightarrows M $ is denoted by $M/ {\cal G} $.

In this paragraph, we show how the quotient spaces in \eqref{isomorphisms} can be characterized in terms of two Lie groupoid orbits. More precisely,  the set of equivalence classes ${PH}_{n}/\sim_{sys}$ (resp. $\Theta_{{CH} _n}/\sim_{\star}$) is the orbit space $PH_n/\mathcal{G}_n$ (resp. $\Theta_{{CH} _n}/\mathcal{H}_n$) of a groupoid $\mathcal{G}_n\rightrightarrows PH_n$ (resp. $\mathcal{H}_n\rightrightarrows\Theta_{{CH} _n}$) which we now construct. The statement \eqref{isomorphisms} is hence equivalent to stating that the orbit spaces $PH_{n}/\sim_{sys}$ and $\Theta_{{CH} _n}/\mathcal{H}_n$ of the two groupoids coincide.

The  total space $\mathcal{G}_n $ of the first groupoid is given by  \begin{multline*}
	\mathcal{G}_n:=\{\left(L,(Q,{\bf B})\right)|L\in GL(2n,\mathbb{R}), (Q,{\bf B})\in PH_n \text{ such that}\\
	\mathbb{J}^TL\mathbb{J}QL^{-1}\text{ is symmetric and positive-definite and } B=\mathbb{J}^TL^T\mathbb{J}L{\bf B}\}.
\end{multline*}
The target and source maps $\alpha,\beta:\mathcal{G}_n\rightarrow PH_n$ are defined by  $\alpha(L,(Q,{\bf B}))$ $:=(\mathbb{J}^TL\mathbb{J}QL^{-1},L{\bf B})$ and $\beta(L,(Q,{\bf B})):=(Q, {\bf B})$.
It can be proved that the orbit space of this groupoid $PH_n/\mathcal{G}_n$ coincides with $PH_n/\sim_{sys}$.

The total space $\mathcal{H}_n $ of the second groupoid is given by
\begin{multline*}
	\mathcal{H}_n:=\big\{\left((P_{\sigma},A),({\bf d},{\bf v})\right)\mid P_{\sigma}\in \mathbb{M}_n \text{ is a permutation matrix}, 
	A\in GL(2n,\mathbb{R}), \\
	({\bf d},{\bf v})\in\Theta_{{CH} _n}, 
	\text{ such that } A^{T}\begin{bmatrix}
		D&0\\0&D
	\end{bmatrix}A{\bf v}=\begin{bmatrix}
		D&0\\0&D
	\end{bmatrix}{\bf v}\\
	\text{ and }
	A\mathbb{J}\begin{bmatrix}
		D&0\\0&D
	\end{bmatrix}=\mathbb{J}\begin{bmatrix}
		D&0\\0&D
	\end{bmatrix}A,
	\text{ with } D={\rm diag}({\bf d})\big\}.
\end{multline*}
The target and source maps $\alpha,\beta:\mathcal{H}_n\rightarrow\Theta_{{CH} _n}$ are defined as $\alpha((P_{\sigma},A),({\bf d},{\bf v}))$ $:=({\bf d},{\bf v})$ and $\beta((P_{\sigma},A),({\bf d},{\bf v})):=(P_{\sigma}{\bf d}, PA{\bf v})$, where $P=\begin{bmatrix}
	P_{\sigma}&0\\0&P_{\sigma}
\end{bmatrix}$.

The orbit space of this groupoid $\Theta_{{CH} _n}/\mathcal{H}_n$ coincides with $\Theta_{{CH} _n}/\sim_{\star}$ and, moreover, the orbit spaces of the Lie groupoids $\mathcal{G}_n\rightrightarrows PH_n$ and $\mathcal{H}_n\rightrightarrows\Theta_{{CH} _n}$ are isomorphic, as desired.

\subsubsection{Acknowledgements} The authors thank Lyudmila Grigoryeva for helpful discussions and remarks and acknowledge partial financial support from the Swiss National Science Foundation (grant number 175801/1) and the School of Physical and Mathematical Sciences of the Nanyang Technological University. DY is funded by the Nanyang President's Graduate Scholarship of Nanyang Technological University.

%
%
%
\bibliographystyle{splncs04}
\bibliography{literature}
%
%
%
%
%
\end{document}